\documentclass{amsart}
\usepackage{amsmath}
\usepackage{amssymb,latexsym}
\usepackage{amsfonts,mathrsfs}
\usepackage[margin=1.8in]{geometry}
\usepackage[all,cmtip]{xy}
\usepackage[colorlinks= true,pdfborder={0 0 0}]{hyperref}
\usepackage{graphicx}
\usepackage{color}
\def\f#1#2{\frac{#1}{#2}}

\def\b{\beta}

\def\({\left (}
\def\){\right )}
\def\<{\left\langle}
\def\>{\right\rangle}

\newtheorem{thm}{Theorem}[section]

\newtheorem{prop}[thm]{Proposition}
\newtheorem{defn}[thm]{Definition}
\newtheorem{rem}[thm]{Remark}

\newcommand{\norm}[1]{\left\Vert#1\right\Vert}
\newcommand{\abs}[1]{\left\vert#1\right\vert}

\newcommand{\eps}{\varepsilon}

\newcommand{\cta}{\theta}

\begin{document}

\title[positive mass on conical manifold]{Positive mass theorem on conical manifold
with small cone angle}
\keywords{ Conical Space; Small Cone Angle; Positive Mass}
\author{Yaoting Gui}
\address{Yaoting Gui, Beijing International Center for Mathematical Research,
Peking University, Beijing, China } \email{ytgui@bicmr.pku.edu.cn}
\thanks {2010 Mathematics Subject Classification. 51F99,31E05}
\thanks{The research was supported by NSFC,  No. 11721101. }

\begin{abstract}
  We prove the positive mass theorem on conical manifold with small cone angle and co-dimensional
  two singularities under the assumption that the ambient manifold admits a spin structure and
  locally conformal flat. 

\end{abstract}

\maketitle
\numberwithin{equation}{section}
\section{\bf Introduction}
\setcounter{equation}{0}
We consider the positive mass theorem on the conical space with small cone angle. On one
hand, the conical singularities with small cone angle shares similar geometric properties
as Riemann surfaces, see \cite{troyanov1991prescribing}. On the other hand, the positive mass
theorem plays a vital role in both physics and mathematics. In the past decades, generalizations
of the the positive mass on singular space also attract many attentions, see \cite{miao2002positive, shi2017scalar} and references therein. Our purpose in this paper 
is to prove some kind of positive mass theorem on the conical space of co-dimensional 2 singularities
with small cone angle under the condition that the ambient manifold be spin and locally conformal flat.
Roughly speaking, we consider a closed Riemannian manifold $M$ with a co-dimensional 2 embedding
submanifold $X$, along which the singularities is conical modelled by
\begin{equation}\label{*}
 X_\beta=\{\mathbb{R}_\beta^2\times\mathbb{R}^{n-2},g_\beta=d\rho^2+(1+\beta)^2\rho^2
    d\theta^2+\sum_{i=3}^{n}d\xi_i^2\}.
\end{equation}
See the precise definition \ref{def2.1}. Basic references about the conical space and positive mass are \cite{bartnik1986mass,
miao2002positive, witten1981new, schoen1981proof, schoen1979proof,parker1982witten, shi2017scalar}
\begin{thm}[Positive Mass Theorem]\label{PMT}
  Let $M$ be a conical manifold with a singular set $X$ of codimension 2 and $\b\in(-1,0)$.
  Suppose the ambient manifold is spin and locally conformal flat with positive Yamabe constant,
  then the mass of $M$ is non-negative. More precisely, let $q\in M\setminus X$ be a regular point.
  Consider the conformal metric $\tilde{g}=G^{\frac{4}{n-2}}g$, where  $G=G_q$ is the Green
  function with respect to the conformal Laplacian $L_\beta=-\f{4(n-1)}{n-2}\Delta_\beta+R_g$. Then
  the new manifold $(\tilde{M}=M\setminus X\cup q,\tilde{g})$ is asymptotically flat of order $n-2$
  with non-negative mass $m(\tilde{g})\geq0$.
\end{thm}
\begin{rem}
The Green function is defined by the following equation,
\begin{equation}\label{Eq34}
L_\beta G=\delta_q,
\end{equation}
where $q\in M\setminus X$ is a regular point.
By the assumption that $M$ is locally conformal flat, we have under a local conformal normal
coordinate, $G$ has the following expansion,
\[
G(x)=dist(q,x)^{2-n}+A+\alpha(x),
\]
where $\alpha(x)=O(r)$ is Lipschitz, and $r$ is the distance with the reference point $q$.
The positive mass theorem asserts that $A$ is nonnegative because $A$ is a positive multiple
of the mass $m(\tilde{g})$. We will prove the positive mass theorem by showing that $A\geq0$.
\end{rem}
\begin{rem}\label{rek36}
  The existence of the Green function and the conformal normal coordinates are guaranteed
  by the local natural. Namely, we can construct the Green function with the pole far away
  from the singularities and then follow the construction presented in \cite{aubin2013some}.
  Moreover, we can choose a smooth point $x_0$ as a reference point and construct the conformal
  normal coordinate as in the smooth case and expand the Green function near $x_0$. For more
  details, see the book \cite{schoen1994lectures} and the survey paper \cite{lee1987yamabe}.
\end{rem}
\begin{rem}
  As illustrated in \cite{ammann2005positive}, the negative Yamabe constant case is easy to settle, so we here additionally assume that the ambient manifold admits positive Yamabe constant. 
\end{rem}
The proof makes use of the Lichnerowicz-Sch\"odinger formula, combined with a construction of the
test spinor. The difficulty is that we need to carefully deal with the estimate near the conical
singularities.
This is indeed follows from the self-adjointness of the Dirac operator on conical manifold.
Note that here we are considering co-dimensional 2 singularities. It is known that the positive
mass is not true any more if the metric is conical along a hypersurface, see \cite{shi2017scalar}.
Hence the co-dimensional 2 singularities is subtle. \medskip

The paper is organized as follows. Section 2 is devoted to some basic materials about conical
manifold and spin structure. And then we move to the proof of our positive mass theorem in section 3.

\textbf{Acknowledge.}
The author would like to express his great thankness to his tutor Prof. Jiayu Li for his kind
and patient guide.
\section{\bf Preliminary}
We first introduce some bases on the conical manifold. The definition is motivated by our study
on the standard model space $X_\b$. Let $(M^n,g_0)$ be a closed smooth manifold equipped with a
Riemannian metric $g_0$ and $\iota: X^{n-2}\hookrightarrow M$ be an $(n-2)$ dimensional embedding
submanifold with the induced metric $g_X=\iota^*g$.
The tubular neighbourhood theorem asserts that there exists an $\eps$-neighbourhood $X_\eps$ inside
its normal bundle $NX$ of $M$. Suppose $M$ admits a metric defined as following. For any $p\in X$,
let $\{x^i\}$ be a normal coordinate system around $p$, and $\{\rho,\theta\}$ a local coordinate in
the normal bundle. Under this coordinate, we may write the metric as following,
$\forall q=(\rho,\cta,x)\in X_\eps$,
\begin{equation}\label{for21}
g=\sigma_{ij}dx_idx_j+d\rho^2+(\b+1)^2\rho^2(d\theta^2+
\sigma_id\theta dx_i+f_{ij}dx_idx_j),
\end{equation}
where $f_{ij}, \sigma_i, \sigma_{ij}$ are all smooth functions of $(\rho,\cta,x)$.
\begin{defn}\label{def2.1}
  We say M is a conical manifold of co-dimensional 2 singularities with small cone angle
  (conical space for sort), if it admits a metric (\ref{for21}) such that the tangent cone
  at a point $p\in X$ is isometric to the standard co-dimensional 2 cone $X_\b$, where
\[
    X_\beta=\{\mathbb{R}_\beta^2\times\mathbb{R}^{n-2},g_\beta=d\rho^2+(1+\beta)^2\rho^2
    d\theta^2+\sum_{i=3}^{n}d\xi_i^2\},
\]
and the total angle of the conical metric less than or equal to $2\pi$. Namely, we assume
$\b\in(-1,0)$ since the total angle of the conical metric is $2\pi(\b+1)$.
\end{defn}
We now recall some basic concepts about spin structure and Green function of the Dirac
operator. See \cite{hijazi2001spectral, friedrich2000dirac, kobayashi1963foundations}.

Suppose $M$ is an oriented Riemannian manifold, probably with boundary. But we do not need
to assume the completeness. In the sequel, we shall apply this construction to the space
$M\setminus X$. Let $SOM$ be the orthonormal frame bundle over the base manifold
$M$. Let $V_\alpha$ be a covering of $M$ and
$\phi_{\alpha\beta}: U_\alpha\cap U_\beta \rightarrow SO(n)$ be the
transition functions. Let $Ad$ be the two-fold universal covering of the spin group $Spin(n)$
to $SO(n)$, then $M$ is called to admit a spin structure if there is a lift
$\eta: SpinM\rightarrow SOM$ such that the diagram commutes
\[
\xymatrix{
                &         SpinM \ar[d]^{\eta}     \\
  V_\alpha\cap V_\beta \ar[ur]^{\tilde{\phi}} \ar[r]_{\phi} & SOM           }
\]
and
\[
Ad\circ\tilde{\phi}_{\alpha\beta}=\phi_{\alpha\beta}\quad and\quad
\tilde{\phi}_{\alpha\beta}\circ\tilde{\phi}_{\beta\gamma}\circ\tilde{\phi}_{\gamma\alpha}=Id.
\]
We say the $Spin(n)$-principal bundle $(SpinM,\eta)$ is a spin structure over $M$.
\begin{defn}
  A complex spinir bundle associated with the spin structure is the complex vector bundle
  \[
  \Sigma M=SpinM\times_\rho \Sigma_n,
  \]
  where $\rho:Spin(n)\rightarrow Aut(\Sigma_n)$ is the complex $Spin(n)$ representation,
  $\Sigma_n$ is the complex spinor space.
\end{defn}
The principal bundle connection on $SOM$ can be naturally lifted to a bundle connection on
$SpinM$, and therefore induces a connection on the associated vector bundle $\Sigma M$, denoted
this connection by $\nabla$. There exists a natural product operation of the tangent bundle on
the associated bundle induced from the Clifford multiplication. Namely, let $s\in\Gamma_U(SOM)$
be a local section defined on an open subset $U\subset M$ and $\tilde{s}$ be its lift on
$\Gamma_U(SpinM)$. Let $\psi=[\tilde{s},\sigma]\in \Gamma(S)$ be a smooth spinor, i.e. smooth
sections of the vector bundle $\Sigma M$. A vector field $X$ is denoted as an equivalent class
$[s,\alpha]$ if we view the tangent bundle as an associated vector bundle of the principal bundle
$SpinM$. Define
\[
cl:TM\otimes \Sigma M\longrightarrow \Sigma M
\]
\[
X\otimes \psi:=[\tilde{s},\alpha]\otimes [\tilde{s},\sigma]\longmapsto
[\tilde{s},\alpha\cdot\sigma]:=X\cdot\psi,
\]
where $"\cdot"$ denote the Clifford multiplication. Since $\Sigma M$ is a complex vector bundle,
there is a natural Hermitian inner product $(,)$ which is preserved by the Clifford multiplication.
That is, for $x\in\Gamma(TM)$ with $|x|=1$, there holds
\[
(x\cdot\phi,x\cdot\psi)=(\phi,\psi),\quad \forall \phi,\psi\in
\Gamma(S),
\]
The associated $L^2(S)$ space is defined
by the completion of smooth spinors $\Gamma(S)$ under the $L^2$-norm, where
\[
||\psi||_{L^2(S)}=\{\psi\in S|\int_M(\psi,\bar{\psi})dv_g<\infty\}.
\]
We are now ready to define the Dirac operator.
\begin{defn}
  The Dirac operator is the composition of the covariant derivative acting on sections of
  $\Sigma M$ with the Clifford multiplication:
  \[
  D=cl\circ\nabla.
  \]
\end{defn}
Locally, choose a local orthonormal frame $e_i$, we have
\[
\xymatrix@C=0.3cm{
 D=cl\circ\nabla: \Gamma(S) \ar[rr]^{\nabla\quad} && \Gamma(T^*M\otimes\Sigma M)
 \ar[rr]^{\quad cl} &&   \Gamma(S) }
\]
by $D\psi=e_i\cdot \nabla_{e_i}\psi$, where $\psi\in S$ is called a spinor field. It is
easy to see that this definition is independent of the choice of the orthonormal bases.
Similarly, we define $D_0$ to be the Dirac operator on smooth spinors with compact support.
\medskip

We list some basic facts about the Dirac operator.
\begin{itemize}\label{dp}
  \item both $D_0$ and $D$ is elliptic;
  \item $D_0$ is symmetric on $L^2(S)$;
  \item $D_0$ has a self-adjoint extension on $L^2(S)$.
\end{itemize}
\begin{rem}
  The self-adjoint extension of $D_0$ is related to the Dirichlet boundary condition if
  $\partial M\neq\emptyset$.
\end{rem}
We first calculate the principal symbol of the Dirac operator. For a locally defined
spinor $\psi=\psi^\alpha\sigma_\alpha$, and an orthonormal basis $e_i$, we have by definition
\[
D\psi=e_i\cdot\nabla_{e_i}(\psi^\alpha\sigma_\alpha)
=\nabla\psi^\alpha\cdot\sigma_\alpha+\psi^\alpha e_i\cdot\nabla_{e_i}\sigma_\alpha,
\]
which implies that
\[
\sigma_P(D)(\xi)\psi=\xi\cdot\psi.
\]
Suppose $\xi\cdot\psi=0$, then multiplied both side by $\xi$ and notice that
$\xi\cdot\xi=-1$, we conclude that $\psi=0$, hence $D$ is elliptic. To see $D_0$ is symmetric,
we choose a parallel orthonormal basis, that is, $\nabla_{e_i}e_j=0$,
then $\forall\psi,\varphi\in C^\infty_0(S)$,
\begin{align}\label{For317}
  (D\psi,\varphi)= & (e_i\cdot\nabla_{e_i}\psi,\varphi)\nonumber \\
   = & -(\nabla_{e_i}\psi,e_i\cdot\varphi) \nonumber \\
   = & -e_i(\psi,e_i\cdot\varphi)+(\psi,e_i\cdot\nabla_{e_i}\varphi)\nonumber\\
   = &  div V+(\psi,D\varphi)
\end{align}
where $V$ is a complex vector field,
and we have used the compatibility between the Hermitian product and the Clifford
multiplication. Integrating both side of (\ref{For317}) and note that $\phi,\psi$ have compact
support, the Stokes formula gives the desired result. The last assertion is obvious by using
the proved facts, see \cite{friedrich2000dirac}.
\begin{rem}
  It is remarkable that there is not always the case to obtain a self-adjoint extension of the
  Dirac operator $D$, see the example \cite{dirac}. The reason is that we may view the conical
  singularities as ideal infinity, and there exist some harmonic spinors which are also $L^2$-
  integrable with non-vanished boundary integral. However, under the assumption that the cone angle
  is small, that is, $\b\in(-1,0)$, we do obtain self-adjoint extension, see \cite{dirac, chou1989criteria}.
  \end{rem}

We shall denote the associated Dirac operator with respect to the conical metric by $\bar{D}$.
Under the above notation, we can also define the Sobolev space $H^1(S)=W^{1,2}(S)$ as the completion
of the smooth spinors under the $H^1$-norm, i.e.
\begin{equation}
\norm{\psi}_{H^1(S)}=\{\psi\in L^2(S)|\int_{M}\abs{\psi}^2+\abs{\bar{D}\psi}^2dv_g\leq\infty\}.
\end{equation}
$W^{1,p}(S)$ can be defined similarly. Let us now derive the so-called Lichnerowicz-Sch\"{o}rdinger
formula, which asserts that
\begin{equation*}\label{ls}
D^2=\nabla^*\nabla+\frac{1}{4}R_g\mathrm{Id}_{\Gamma(S)},
\end{equation*}
where $\nabla^*$ is the formal adjoint and $R_g$ is the scaler curvature. To see this, we
choose a local orthonormal basis $\{e_i\}$ satisfying $\nabla_{e_i}e_j=0$, then we have
\begin{align*}
  D^2 = & (e_i\cdot\nabla_{e_i})(e_j\cdot\nabla_{e_j})=e_i\cdot\nabla_{e_i}e_j\cdot\nabla_{e_j}
  +e_i\cdot e_j\cdot\nabla_{e_i}\nabla_{e_j} \\
  = & -\nabla_{e_i}\nabla_{e_i}+\Sigma_{i<j}e_i\cdot
  e_j(\nabla_{e_i}\nabla_{e_j}-\nabla_{e_j}\nabla_{e_i}) \\
  = & -\nabla_{e_i}\nabla_{e_i}+\Sigma_{i<j}e_i\cdot
  e_j\cdot R_{e_i,e_j} \\
  = & -\nabla_{e_i}\nabla_{e_i}+\frac{1}{2}\Sigma_{i,j}e_i\cdot
  e_j \cdot R_{e_i,e_j}.
\end{align*}
We claim that
\[
\Sigma_{i,j}e_i\cdot
  e_j \cdot R_{e_i,e_j}=\frac{1}{4}R_g\rm{Id}_{\Gamma(S)}
  \]
and we further note that
\[
(\nabla_{e_i}\nabla_{e_i}\psi,\varphi)=e_i(\nabla_{e_i}\psi,\varphi)
-(\nabla_{e_i}\psi,\nabla_{e_i}\varphi)
\]
Since $e_i(\nabla_{e_i}\psi,\varphi)$ is a divergence of some complex vector field, we conclude
that
\[
-\int_M(\nabla_{e_i}\nabla_{e_i}\psi,\varphi)=\int_M(\nabla_{e_i}\psi,\nabla_{e_i}\varphi),
\]
or
\[
-\int_M(\nabla^*\nabla\psi,\varphi)=\int_M(\nabla\psi,\nabla\varphi).
\]
$\forall\phi\in\Gamma(S)$.
We then obtain the desired result. in the conical setting, we have a similar equation (\ref{ls})
\begin{equation*}
  D^2=-\nabla^*\nabla+\frac{1}{4}R_g\mathrm{Id}_{\Gamma(S)},\quad \mbox{in}
  \quad M\setminus X,
  \end{equation*}
 where all the operators are taken with respect to the conical metric. In the sequel, we only deal
 with the conical metric, for simplicity, we will still denote these operators without notational
 change. 

\section{\bf Proof of the Positive mass}
In this section, we prove the positive mass theorem on the conical space with small cone angle and
admits a spin structure. We will construct the Green function for the conical Dirac operator $D$.
This is an analogy of the construction of Dirac operator on the smooth closed manifold, see \cite{ammann2005positive}. For simplicity, we will assume that $M$ is locally conformal flat, that
is, $\forall p\in M$, there exists a neighbourhood $U_p$ such that there is a flat metric in the
conformal class $[g]$.
\begin{prop}\label{prop319}
  Let $x_0\in M\setminus X$ be a regular point, denote $\{x^i\}$ the local coordinate in
  a small neighbourhood $B_\delta(x_0)$ with $\delta<\min\{i_{x_0},d(x_0,X)\}$, where $i_{x_0}$ is
  the injectivity radius of the point $x_0$ and further
  $\delta$ is so chosen that the spinor  bundle over $B_\delta(x_0)$ is trivial. Let $\psi_0$
  be a constnat spinor, then there exists a $D$-harmonic spinor $\psi$ on $M\setminus X\cup
  \{x_0\}$ such that
  \begin{equation}\label{for31}
  \psi|_{B_\delta(x_0)}=\frac{x}{r^n}\cdot\psi_0+\theta(x),
  \end{equation}
  where $\theta(x)$ is a smooth spinor on $B_\delta(x_0)$
\end{prop}
\begin{proof}
  Consider a cut-off function $\zeta$,
  \[
  \zeta=1\quad in\quad B_\frac{\delta}{2}(x_0),\quad \zeta=0\quad on\quad M\setminus B_\delta.
  \]
  Define $\Phi=\frac{x}{r^n}\zeta(x)\cdot\psi_0$, where $"\cdot"$ denote the Clifford
  multiplication. Since on $B_\delta$ the metric is flat, $\bar{D}=e_i\cdot \partial_i$, it is
  straightforward to check that $D(\frac{x}{r^n}\cdot\psi_0)=0$. Therefore
  $\Phi$ is $D$-harmonic near $x_0$. Set $f=D\Phi$, then
  $f$ vanishes near $x_0$ and hence bounded, which can be
extended to become a bounded smooth spinor on $M\setminus X$ and is $L^2$-intergrable. As
$D$ is self-adjoint, we can solve the equation $D\theta=-f$ in $M\setminus X$.
\begin{rem}\label{rem3.2}
  The self-adjointness of the Dirac operator implies that a general smooth spinor $\psi$ which
  is $L^2$-integrable together with its derivative $D\psi$ can be approximated by smooth spinors
  with compact support in $M\setminus X$. Indeed, we always have $Dom(D_0)\subset Dom(\bar{D})
  \subset Dom(D_0^*)$, where for simplicity, we use the same notation $D_0$ to also denote its
  $L^2$-closure. Now the self-adjointness of $D$ implies that $Dom(D_0)= Dom(\bar{D})= Dom(D_0^*)$,
  which shows that a general smooth spinor can ba approximated by smooth spinors with compact support.
\end{rem}
\end{proof}
\begin{proof}[Proof of Positive Mass Theorem \ref{PMT}]
Consider the conformal metric $\tilde{g}=G^{\frac{4}{n-2}}g$, In the following, all
the notation is taken with respect to the metric $\tilde{g}$. By the formula of conformal
change for the scaler curvature, we have $R_{\tilde{g}}=0$. And furthermore, if $\tilde{g}=e^{2u}g$,
\[
\tilde{D}(e^{-\frac{n-1}{2}u}\tilde{\psi})=e^{-\frac{n+1}{2}u}\widetilde{D\psi},
\]
thus $\tilde{\psi}=G^{-\frac{n-1}{n-2}}\psi$ is $\tilde{D}$-harmonic. We have by the
Sch\"{o}rdinger-Lichnerowicz formula,
\[
0=\tilde{D}^2\tilde{\psi}=\tilde{\nabla}^*\tilde{\nabla}\tilde{\psi}
+\frac{1}{4}R_{\tilde{g}}\tilde{\psi}\quad in\quad M\setminus X.
\]
Let $B_\eps(q)$ be a small ball around $q$ such that $B_\eps(q)\cap X=\emptyset$
and $U_\eps$ be a small neighbourhood of $X$.
In particular, we choose the test spinor $\psi$ as in (\ref{prop319}) with the pole
$x_0\in B_\epsilon(q)$ and integrate over $M\setminus B_\eps(q)\cup U_\eps$,
\[
0=\int_{M\setminus ({B_\eps(q)\cup U_\epsilon})}(\tilde{\psi},\tilde{\nabla}^*
\tilde{\nabla}\tilde{\psi})dv_{\tilde{g}}.
\]
Applying the stokes formula, we find that
\[
\int_{M\setminus (B_\eps(q)\cup U_\eps)}|\nabla\tilde{\psi}|^2
=\int_{\partial B_\eps(q)}(\partial_{\tilde{\nu}}\tilde{\psi},\tilde{\psi})
+\int_{\partial U_\eps}(\partial_{\tilde{\nu}}\tilde{\psi},\tilde{\psi}).
\]
It follows that
\begin{equation}\label{In329}
\int_{\partial B_\eps(q)}\partial_{\tilde{\nu}}(\tilde{\psi},\tilde{\psi})
+\int_{\partial U_\eps}\partial_{\tilde{\nu}}(\tilde{\psi},\tilde{\psi})\geq0.
\end{equation}
A simple calculation gives that
\begin{equation*}\label{}
\int_{\partial B_\eps(q)}\partial_{\tilde{\nu}}(\tilde{\psi},\tilde{\psi})
=2(n-1)\omega_{n-1}A+O(\eps).
\end{equation*}
Indeed,
\[
\tilde{\nu}=G^{-\frac{2}{n-2}}\nu=G^{-\frac{2}{n-2}}\frac{\partial}{\partial r},
\]
\[
d\sigma_{\tilde{g}}=G^{\frac{2(n-1)}{n-2}}d\sigma_g=G^{\frac{2(n-1)}{n-2}}
\epsilon^{n-1}d\sigma \quad on\quad \partial B_\epsilon(q),
\]
here $d\sigma$ is the usual spherical measure on $S^{n-1}$,
because the metric $g$ is flat, and we have chosen sufficiently small $\epsilon$
so that $\partial B_\epsilon(q)$ is smooth. We may assume $|\psi_0|=1$, note that,
\begin{align*}
  |\tilde{\psi}|^2_{\tilde{g}} & = G^{-\frac{2(n-1)}{n-2}}|\psi|^2_g \\
   & = (r^{2-n}+A+\alpha)^{-\frac{2(n-1)}{n-2}}|\frac{x}{r^n}
  \cdot\psi_0+\theta|^2_g \\
   & = (1+Ar^{n-2}+\alpha r^{n-2})^{-\frac{2(n-1)}{n-2}}(1+2r^{n-1}
   Re(\frac{x}{r}\cdot\psi_0,\bar{\theta})+r^{2n-2}|\theta|^2_g).
\end{align*}
We compute as following,
\begin{align*}
  \frac{\partial}{\partial r}|\tilde{\psi}|^2 & =-2(n-1)\Big(1+Ar^{n-2}
  +\alpha r^{n-2}\Big)^{-\frac{3n-4}{n-2}}\Big(Ar^{n-3}+O(r^{n-2})\Big) \\
   &  \quad\Big(1+2r^{n-1}Re(\frac{x}{r}\cdot\psi_0,\bar{\theta})+r^{2n-2}|\theta|^2\Big)
    +\Big(1+Ar^{n-2}+O(r^{n-1})\Big)^{-\frac{2(n-1)}{n-2}}\\
   & \quad \Big((2n-2)r^{n-2}Re(\frac{x}
   {r}\cdot\psi_0,\bar{\theta})+(2n-2)r^{2n-3}|\theta|^2\Big).
\end{align*}
It follows that
\[
\frac{\partial}{\partial r}|\tilde{\psi}|^2=(1+Ar^{n-2}+O(r^{n-1}))^{-\frac{2(n-1)}{n-2}}
\bigg(\frac{-2(n-1)Ar^{n-3}+O(r^{n-2})}{1+Ar^{n-2}+O(r^{n-1})}+O(r^{n-2})\bigg),
\]
hence
\[
\partial_{\tilde{\nu}}|\tilde{\psi}|^2=G^{-\frac{2}{n-2}}\Big(2(n-1)Ar^{n-3}+O(r^{n-2})\Big).
\]
We then obtain
\begin{align*}
  0 & \leq \int_{\partial B_\epsilon(q)}\partial_{\tilde{\nu}}
  |\tilde{\psi}|^2d\sigma_{\tilde{g}}\\
   & = \int_{S^{n-1}}2(n-1)G^2\epsilon^{n-1}(A\epsilon^{n-3}+O(\epsilon^{n-2})) \\
   & = \int_{S^{n-1}}2(n-1)(\epsilon^{2-n}+A+O(\epsilon))^2\epsilon^{n-1}
   (A\epsilon^{n-3}+O(\epsilon^{n-2})) \\
   & = \int_{\partial B_1}2(n-1)(\epsilon^{3-n}+A^2+O(\epsilon))
   (A\epsilon^{n-3}+O(\epsilon^{n-2})) \\
   & = 2(n-1)\omega_{n-1}A+O(\epsilon).
\end{align*}
It suffices to prove that the second integral in (\ref{In329}) vanishes. By Remark \ref{rem3.2}, we may
first assume that the spinor $\cta$ has compact support in $M\setminus X$, which immediately
implies that the second integral in (\ref{In329}) vanishes. The general case is a simple approximation
by smooth spinors with compact support. We now complete the proof.
\end{proof}
\begin{rem}
  The proof can be indeed generalized to more general space, namely pseudo-manifold in the sense of \cite{chou1989criteria}. This is because the essential point that we need is the self-adjointness of
  the Dirac operator, which is true in such pseudo-manifold.
\end{rem}

\bibliography{ref}
\bibliographystyle{plain}

\end{document}